\input amstex
\documentstyle{amsppt}
\magnification 1200
\NoRunningHeads

\topmatter
\title 
On limit theorems for fields of martingale differences
\endtitle
\author 
Dalibor Voln\'y 
\endauthor
\affil
Laboratoire de Math\'ematiques Rapha\"el Salem, UMR 6085, Universit\'e de Rouen - Normandie, France
\endaffil
\thanks
I greatly thank the two unknown referees for detailed and careful reading of the manuscript and for many helpful comments and remarks.
\endthanks

\abstract
We prove a central limit theorem for stationary multiple (random) fields of martingale differences $f\circ T_{\underline{i}}$, 
$\underline{i}\in \Bbb Z^d$, where $T_{\underline{i}}$ is a $\Bbb Z^d$ action.  In most cases the  multiple (random) fields of martingale differences
 is given by a completely  commuting filtration.
A central limit theorem proving convergence to a normal law has been known for Bernoulli random fields and in \cite{V15} this result was extended 
to random fields where one of generating transformations is ergodic.

In the present paper it is proved that a convergence takes place  always and the limit law is a mixture of normal laws. If the  $\Bbb Z^d$ action
is ergodic and $d\geq 2$, the limit law need not be normal.

For proving the result mentioned above, a generalisation of McLeish's CLT for arrays $(X_{n,i})$ of
martingale differences is used. More precisely, sufficient conditions for a CLT are found in the case when the sums $\sum_i X_{n,i}^2$ converge only
in distribution. 

The CLT is followed by a weak invariance principle. It is shown that central limit theorems and invariance principles using martingale approximation
remain valid in the non-ergodic case.
\endabstract

\endtopmatter

\subheading{1. Introduction and central limit theorems}
In the study of limit theorems for dependent random variables an important role has been played by the central limit theorem for 
ergodic sequences of matingale differences, found independently by P\. Billingsley and I.A\. Ibragimov (cf\. \cite{B61}, \cite{I}). In the non-ergodic case
 the theorem
remains true but the limit law is a mixture of non normal laws (cf\.  \cite{HaHe}), the result has been proved and reproved many times; in \cite{V89}
ergodic decomposition of an invariant measure, as in this paper, is used. 
Here, we will study the CLT for random fields of martingale differences.
By a random field we 
understand a field of random variables $f\circ T_{\underline{i}}$, $\underline{i} \in \Bbb Z^d$, on a probability space  $(\Omega, \Cal A, \mu)$ 
where $f$ is a measurable function on $\Omega$ and $ T_{\underline{i}}$ are automorphisms of $(\Omega, \Cal A, \mu)$ for which
$  T_{\underline{i}} \circ  T_{\underline{j}} =  T_{\underline{i+j}}$ ($( T_{\underline{i}})$ is a $\Bbb Z^d$ action). We denote $f\circ T_{\underline{i}}
= U_{\underline{i}}f$. Recall that the $\Bbb Z^d$ action 
is {\it ergodic} if the only measurable sets $A$ for which $ T_{\underline{i}}A = A$ are of measure zero or one (for $d=1$ we speak of an
ergodic transformation $T = T_1$). 
Limit theorems for random fields of martingale differences have been studied by e.g\. Basu and Dorea \cite{BaDo}, or by Nahapetian \cite{N}.
An approach using so called projection method was used in \cite{D}. 
For a random field, there are several non-equivalent definitions of martingale differences.

For martingale approximations and for weak invarince principles
we will use the notion of orthomartingales as in e.g\. \cite{Go09}, \cite{WaW} (cf\. also \cite{K}):

For $\underline{i} =(i_1,\dots,i_d)\in \Bbb Z^d$, $\underline{i} \leq \underline{j}$  means that $i_k\leq j_k$, $k=1,\dots, d$;  
$\underline{i}\wedge \underline{j} = (\min\{i_1, j_1\}, \dots, \min\{i_d, j_d\})$.  $(\Cal F_{\underline{i}})_{\underline{i}\in\Bbb Z^d}$ is a  
{\it completely commuting invariant filtration} if
\roster
\item"(i)" $\Cal F_{\underline{i}}= T_{-\underline{i}} \Cal F_{\underline{0}}$ for all $\underline{i} \in \Bbb Z^d$,
\item"(ii)" $\Cal F_{\underline{i}} \subset F_{\underline{j}}$ for $\underline{i} \leq \underline{j}$,
\item"(iii)" $\Cal F_{\underline{i}} \cap \Cal F_{\underline{j}} = \Cal F_{\underline{i}\wedge 
\underline{j}}$, $\underline{i}, \underline{j}\in \Bbb Z^d$, 
\item"(iv)"  $E\Big( E(f | \Cal F_{\underline{i}}) \big| \Cal F_{\underline{j}}\Big) = E(f | \Cal F_{\underline{i}\wedge \underline{j}})$,
for every integrable function $f$.
\endroster
By $\Cal F_l^{(q)}$, $1\leq q\leq d$, $l\in \Bbb Z$, we denote the $\sigma$-algebra generated by the union of all $\Cal F_{\underline{i}}$ 
with $i_q\leq l$.
For $d=2$, $\Cal F_{\infty, j} = \Cal F_j^{(2)}$ denotes the $\sigma$-algebra generated by the union of all $\Cal F_{i, j}$, 
$i\in \Bbb Z$, and in the same way we define $\Cal F_{i, \infty}$. \newline
By $P_l^{(q)}$, $1\leq q\leq d$, $l\in \Bbb Z$, we denote the operator in $L^p$, $1\leq p<\infty$, which sends $f\in L^p$ to
$E(f | \Cal F_l^{(q)}) - E(f | \Cal F_{l-1}^{(q)})$. Notice that for $p=2$,  $P_l^{(q)}$ is the orthogonal projection onto the Hilbert space
$L^2(\Cal F_l^{(q)}) \ominus L^2(\Cal F_{l-1}^{(q)})$.  $P_l^{(q)}$ are mutually commuting and
idempotent operators. For $\underline{i} = (i_1,\dots,i_d)$ we define $P_{\underline{i}} = \Pi_{q=1}^d P_{i_q}^{(q)}$. 
In $L^2$, $P_{\underline{i}}$ is the orthogonal projection
onto $\underset 1\leq q\leq d \to{\bigcap} L^2(\Cal F_{i_q}^{(q)}) \ominus L^2(\Cal F_{i_q-1}^{(q)})$.  The functions  $P_{\underline{i}}f$ are called 
{\it martingale differences}. 
Let $\underline{e_i}$ be the vector from $\Bbb Z^d$ with 1 at $i$-th coordinate and 0 elsewhere.
It can be noticed that if $f= P_{\underline{i}}f$ then for $1\leq q\leq d$, $(f\circ T_{\underline{e_q}}^j)_j$ is a martingale difference sequence with 
respect to the filtration $(\Cal F_j^{(q)})_j$. For proofs and for more properties of the operators  $P_{\underline{i}}$ the reader can consult 
 (e.g.) \cite{WaW}, \cite{VWa}.

By central limit theorem we understand weak convergence of the distribution of
$(1/\sqrt{n_1\dots n_d}) \sum_{i_1=1}^{n_1} \dots \sum_{i_d=1}^{n_d} f\circ T_{(i_1,\dots,i_d)}$ as $\min\{n_1,\dots,n_d\} \to\infty$.

 In one dimensional case, the central limit theorem for stationary martingale differences led to many results using
martingale approximations. The same holds true for higher dimension. Pioneering results were given in \cite{WaW} where under a multiparameter version 
of reinforced Maxwell-Woodroofe 
condition  a CLT and weak invariance principle were proved, and in \cite{Go09} where the martingale-coboundary representation was studied 
(an interesting application of Gordin's result can be found in \cite{DeGo}). The results from  \cite{WaW} were improved in \cite{VWa} where the Hannan's condition was generalised 
to random fields, in \cite{PZ} where the CLT was proved under  ``classical"  Maxwell-Woodroofe condition and in \cite{Gi} where the weak invariance
principle (under the same assumptions) was proved. The martingale-coboundary representation has been studied (after \cite{Go09}) in \cite{EGi}, 
\cite{V16}, and \cite{Gi}. More results were published in  (e.g\.) \cite{BiDu}.

For orthomartingale differences, ergodicity of the $\Bbb Z^d$ action (with $d\geq 2$) does not guarantee a convergence to a normal law.
An example can be found in the article by Wang and Woodroofe \cite{WaW}; we use its presentation from \cite{V15}. The idea of the example 
has probably appeared already in the Ph.D thesis of Hillel Furstenberg.
\medskip

\underbar{Example.}  Let the probability space $(\Omega, \Cal A, \mu)$ be a product of  $(\Omega_1, \Cal A_1, \mu_1)$ and  
$(\Omega_2, \Cal A_2, \mu_2)$. On  $(\Omega_i, \Cal A_i, \mu_i)$ there is a bimeasurable and measure preserving bijection $T_i$, $i=1,2$.
On  $(\Omega, \Cal A, \mu)$ we define an action of $\Bbb Z^2$ by $T_{i,j}(x,y) = (T_1^ix, T_2^jy)$. Let us suppose that there is a random variable 
$e_i$ such that $e_i\circ T_i^j$, $j\in \Bbb Z$ are iid $\Cal N(0,1)$ random variables generating the $\sigma$-algebra $\Cal A_i$, $i=1,2$.
$T_{i,j}$ is then an ergodic $\Bbb Z^2$ action (cf\. \cite{V15}). For $e =(e_1, e_2)$, $(e\circ T_{i,j})_{(i,j)\in \Bbb Z^2}$ is then a field of
martingale differences for the natural filtration (cf\. \cite{V15}). As we can easily see, for any integers $n,m \geq 1$, the sum
$(1/\sqrt{nm}) \sum_{i=1}^n \sum_{j=1}^m e\circ T_{i,j}$ is distributed as a product $XY$ of two independent $\Cal N(0,1)$ random variables.
\medskip

A convergence to a normal law is guaranteed if the $\Bbb Z^d$ action is Bernoulli, i.e\. the $\sigma$-algebra $\Cal A$ is generated by
iid random variables $e\circ T_{\underline{i}}$,  $\underline{i}\in \Bbb Z^d$. This assumption was used e.g\. in \cite{WaW} or \cite{BiDu}.

In  \cite{V15} it has been proved that if one of the transformations $T_{\underline{e_i}}$ (recall that $\underline{e_i}$ is the vector from
$\Bbb Z^d$ with 1 at $i$-th coordinate and 0 elsewhere) is ergodic, then for a  random field of square integrable martingale differences
the central limit theorem takes place with a normal law for limit. In \cite{CDV}, the result was extended to reversed martingales and an 
invariance principle was proved.

The problem we study in this paper is whether without any ergodicity 
assumption there still is a convergence to a limit law (as in the example above), or whether it can happen that there is no convergence at all. Theorem 1 gives
a positive answer showing that there always is a limit law which is a mixture of normal laws.

In the theorem we will use a notion of a field of martingale differences which is weaker than the notion defined above. 

Let $T_{\underline{i}}$ be a $\Bbb Z^d$ action on $(\Omega, \Cal A, \mu)$. We say that $(\Cal F_i^{(q)})_{i\in\Bbb Z,\, 1\leq q\leq d}$ is a 
{\it multiple filtration} if for every $1\leq q\leq d$, 
$(\Cal F_i^{(q)})_{i\in\Bbb Z}$ is a filtration with 
$$ \gather
  \Cal F_i^{(q)}  \subset T_{-\underline{e_q}} \Cal F_i^{(q)} =  \Cal F_{i+1}^{(q)},\,\, i\in\Bbb Z,\,\, 1\leq q\leq d, \\
  T_{\underline{e_{q'}}} \Cal F_i^{(q)} = \Cal F_i^{(q)}\,\,\,\text{ for all}\,\,\, 1\leq q'\leq d, \,\,q'\neq q,\,\, i\in\Bbb Z.
  \endgather
$$
 If $f\in L^1$ and there is a multiple filtration $(\Cal F_i^{(q)})_{i\in\Bbb Z,\, 1\leq q\leq d}$ such that 

$f=E(f\,|\, \Cal F_0^{(q)}) - E(f\,|\, \Cal F_{-1}^{(q)})$ for all $1\leq q\leq d$, \newline
we say that  $(f\circ T_{\underline{i}})_{\underline{i}}$
is a {\it multiple field of martingale differences}. If $f\in L^2$ we then have $f\in L^2(\Cal F_0^{(q)}) \ominus  L^2(\Cal F_{-1}^{(q)}) =
U_{e_{q'}} \Big( L^2(\Cal F_0^{(q)}) \ominus  L^2(\Cal F_{-1}^{(q)}) \Big)$ for $1\leq q'\leq d$, $q'\neq q$.

\proclaim{Theorem 1} Let $f\in L^2$  be such that $(f\circ T_{\underline{i}})_{\underline{i}}$ is a multiple field of martingale differences.
If $n_j \to\infty$, $j=1,\dots, d$ then the random variables
$$
  \frac1{\sqrt{n_1\dots n_d}} \sum_{i_1=1}^{n_1} \dots \sum_{i_d=1}^{n_d} f\circ T_{(i_1,\dots,i_d)}
$$
converge in distribution to a law with characteristic function 
$E\exp(-\eta^2t^2/2)$ for a positive random variable $\eta^2$ such that $E\eta^2 = \|f\|_2^2$. 
The random variables
$$
  \frac1{n_1\dots n_d} \sum_{i_1=1}^{n_1} \Big(  \sum_{i_2=1}^{n_2} \dots \sum_{i_d=1}^{n_d} f\circ T_{(i_1,\dots,i_d)}\Big)^2
$$ 
converge in distribution to $\eta^2$.
\endproclaim
\medskip

In order to prove Theorem 1 we use a version of McLeish's theorem \cite{Mc}. McLeish's theorem generalises the Billingsley-Ibragimov's theorem
to triangular arrays $(X_{n,i})_{1\leq i\leq k_n}$ of martingale differences which need not be stationary. Under the assumptions (i), (ii) from the Theorem 2 
below on $\max_i |X_{n,i}|$ and the assumption of convergence in probability of $\sum_{i=1}^{k_n} X_{n,i}^2$ to a constant it gives a CLT for the sums $\sum_{i=1}^{k_n} X_{n,i}$. The book \cite{HaHe} brings several new versions of the theorem, in particular it is  generalised to the case when the 
filtrations $(\Cal F_{n,i})_i$ of the sequences $( X_{n,i})_i$ are nested and $\sum_{i=1}^{k_n} X_{n,i}^2$ converges to a random variable $\eta^2$ 
which is measurable with respect to the intersection of all $\Cal F_{n,i}$. In  \cite{GHu} the conditions (i), (ii) were replaced by 
$\max_{1\leq i\leq k_n} |X_{n,i}| \to 0$ in $L^2$ and in  \cite{L} the $L^2$ convergence was replaced by $L^1$ convergence.
Another generalisation/version of McLeish's theorem was given in \cite{PZ}.
The convergence, nevertheless, is always in probability and a counterexample in \cite{HaHe} shows that a convergence of $\sum_{i=1}^{k_n} X_{n,i}^2$
in distribution is not sufficient without strenghtening of other assumptions. Our Theorem 2 and its corollary, Proposition 3, bring a new version of McLeish's theorem where  the sums $\sum_{i=1}^{k_n} X_{n,i}^2$ converge in distribution only.

\proclaim{Theorem 2} Let  $X_{n,j}$, $j=1,\dots,k_n$, be an array of martingale 
differences with respect to increasing filtrations $(\Cal F_j^n)_{j\geq 0}$, $n=1,2, \dots$,  such that
\roster
\item"(i)" $\max_{1\leq i\leq k_n} |X_{n,i}| \to 0$ in probability,
\item"(ii)" there is an $L<\infty$ such that $E \max_{1\leq j\leq k_n} X_{n,j}^2 \leq L$ for all $n$,
\item"(iii)" there exist $1\leq \ell(n) \leq k_n$ and $\Cal F_{\ell(n)}^n$-measurable random variables $\eta_n^2$  such that
$\sum_{j=1}^{k_n} X_{n,j}^2 - \eta_n^2\to 0$ in probability,
\item"(iv)" there exists a random variable $\eta^2$ such that $ \eta_n^2\to \eta^2$  in distribution,
\item"(v)" $E(T_n(t) | \Cal F_{\ell(n)}^n) \to 1$ in $L^1$ for every $t\in \Bbb R$, where \newline
$T_n(t) = \prod_{j=1}^{k_n} (1 + itX_{n,j}).$
\endroster
Then the sums $\sum_{j=1}^{k_n} X_{n,j}$ converge in distribution to a law with characteristic function 
$E\exp(-\eta^2t^2/2)$.
\endproclaim

As a corollary we get the next proposition.

\proclaim{Proposition 3} Let  $X_{n,j}$, $j=1,\dots,k_n$, be an array of martingale 
differences with respect to increasing filtrations $(\Cal F_j^n)_{j\geq 0}$, $n=1,2, \dots$,  such that assumptions
(i) - (iv) are satisfied and 
\roster
\item"(vi)" the random variables  $\eta_n^2$  are $\Cal F_0^n$-measurable \newline
or
\item"(vii)" the sequences $(X_{n,j})_j$ are strictly stationary and $\eta_n^2$ are measurable with respect to the $\sigma$-algebras of invariant sets.
\endroster
Then the conclusion of Theorem 2 holds.
\endproclaim

\subheading{2. Proofs of Theorems 1,2 and of Proposition 3}

\demo{Proof  of  Theorem 2} As in \cite{HaHe} (cf\. the proof of Theorem 3.2) we for a $C>0$ define the stopping time
$$
  J_n = J_{C,n} = \cases j \,\,&\text{if}\,\, 1\leq j\leq k_n,\,\, \sum_{u=1}^{j-1} X_{n,u}^2 \leq C, \,\, \sum_{u=1}^j  X_{n,u}^2 >C,\\
                                         k_n\,\,&\text{if}\,\, \sum_{u=1}^{k_n} X_{n,u}^2 \leq C
  \endcases
$$
and replace $X_{n,j}$ by martingale differences
$$
  X'_{n,j} = X_{n,j} I_{j\leq J_n}, \,\,\,\, j=1,\dots,k_n.
$$
Notice that  $j\leq J_n$ iff $\sum_{u=1}^{j-1} X_{n,u}^2 \leq C$. We denote
$$
  \eta_C^2 = \eta^21_{\eta^2 \leq C} + C1_{\eta^2 >C}, \quad \eta_{n,C}^2 = \eta_n^21_{\eta_n^2 \leq C} + C1_{\eta_n^2 >C}, \quad n\geq 1.
$$
 From (iv) it follows
$ \eta_{n,C}^2\to \eta_C^2$ in distribution;
from (i), (iii) it follows that 
$\sum_{j=1}^{k_n} {X'_{n,j}}^2 - \eta_{n,C}^2 \to 0$ in probability.

We have verified (iii), (iv) for $ X'_{n,j}$ and  $\eta_C^2$, $\eta_{n,C}^2$. (i) and (ii) follow from $|X'_{n,j}| \leq |X_{n,j}|$. \newline
To see (v), notice that $J_n$ is a finite stopping time; by $\Cal F_{J_n}^n$ we denote the corresponding $\sigma$-algebra. Then 
$\|E(T'_n(t) | \Cal F_{\ell(n)}^n) -1\|_1 = \|E\big(E(T_n(t) | \Cal F_{\ell(n)}^n) -1\,|\,\Cal F_{J_n}^n\big) \|_1$ and we get (v) for $ X'_{n,j}$ by the 
contraction property of conditional expectation.

From (iii), (iv) we deduce
$$
  \lim_{C\to\infty} \limsup_{n\to\infty} \mu(\exists \,1\leq j\leq k_n,\,\, X'_{n,j} \neq X_{n,j}) = 0.
$$
It is thus sufficient to prove the theorem for $X'_{n,j}$, $\eta_C^2$ and $\eta_{n,C}^2$. For simplicity of notation we denote
$ X'_{n,j} = X_{n,j}$ for all $n,j$ and $\eta^2 = \eta_C^2$, $\eta_n^2 = \eta_{n,C}^2$. 
Using (ii) we then get
$$
  E |T_n(t)|^2 = E \prod_{j=1}^{k_n} |1 + itX_{n,j}|^2 = 
   E \prod_{j=1}^{k_n} (1 + t^2X_{n,j}^2)
  \leq e^{t^2C} (1+t^2L).   \tag1
$$
\medskip

As in \cite{Mc} we use the equality 
$$
  e^{ix} = (1+ix)\exp(-\frac12 x^2 +r(x))
$$
where $|r(x)| \leq |x|^3$ for $|x| \leq 1$ ($x\in \Bbb R$) and define
$$
  U_n(t) = \exp\Big(-\frac12 t^2 \sum_{j=1}^{k_n} X_{n,j}^2 + \sum_{j=1}^{k_n} r(t X_{n,j})\Big), \,\,\,t\in \Bbb R.
$$
For $I_n(t) = \exp(itS_n) = T_n(t) U_n(t)$ we then have
$$
  I_n(t) = T_n(t)\Big(U_n(t) - \exp(-\frac12 t^2\eta_n^2)\Big) + 
    T_n(t)\exp(-\frac12 t^2\eta_n^2).
$$
In the same way as in \cite{Mc} we from (i), (iii) deduce that $ \sum_{j=1}^{k_n} r(t X_{n,j}) \to 0$ in probability. Therefore,
$U_n(t) - \exp(-\frac12 t^2\eta_n^2) \to 0$ in probability. $I_n(t)$ are uniformly bounded, hence uniformly integrable; by (1), $T_n(t)$ are
uniformly untegrable. Therefore,
$$
  T_n(t)\Big(U_n(t) - \exp(-\frac12 t^2\eta_n^2)\Big) \to 0\,\,\,\text{in}\,\, L^1.
$$
Because  $\eta_n^2$ are $\Cal F_{\ell(n)}^n$-measurable (cf\. (iii)) we have
$$
  E \exp(-\frac12 t^2\eta_n^2) [T_n(t) - 1] = E  \Big[\exp(-\frac12 t^2\eta_n^2) \big[ E\big(T_n(t) \,|\, \Cal F_{\ell(n)}^n\big) - 1 \big] \Big];
$$
by (v), uniform integrability of $T_n(t)$, and (iv) we deduce
$$
   E T_n(t) \exp(-\frac12 t^2\eta_n^2) \to E \exp(-\frac12 t^2\eta^2)
$$
which finishes the proof.
\enddemo
\qed

\demo{Proof  of  Proposition 3} If the assumptions (i) - (iv) and (vi) are fulfilled, we for all $n$ define $\ell(n)=0$; in the same way as in
\cite{Mc} or \cite{HaHe} we can see that $E(T_n(t)\,|\, \Cal F_0^n) = 1$ for all $n, t$. 

 Let the assumptions (i) - (iv) and (vii) be fulfilled. For each of the sequences $(X_{n,j})_j$ there exists a measure preserving transformation
$T=T_{(n)}$ such that $X_{n,j} = X_{n,0}\circ T^j$ for all $j$. By $\Cal I = \Cal I_n$ we denote the $\sigma$-algebra of $T$-invariant
(measurable) sets. As shown in \cite{V87}, for every $\Cal F_\infty^n$-measurable and integrable function $f$ we have
$E(f\,|\, \Cal F_k^n) = E(f\,|\, \Cal F_k^n \vee \Cal I)$ where $\Cal F_k^n \vee \Cal I$ is the $\sigma$-algebra generated by $\Cal F_k^n \cup \Cal I$.
We thus can replace $\Cal F_k^n$ by $\Cal F_k^n \vee \Cal I$ and hence get (vi).

In both cases the result follows from Theorem 2.
\enddemo
\qed

\demo{Proof of  Theorem 1}

First, we prove the theorem for $d=2$. The general case can  be  proved by induction. 
Recall that we denote $U_{i,j} f = f\circ T_{i,j}$.
\medskip

 If $f$ is not bounded, for a $K>0$ we define 
$$
  f' = P_{0,0} (f1_{|f|\leq K}),\quad f'' = f-f'.
$$
Because $U_{i,j}f$ and  $U_{i,j}f'$ are martingale differences,  $U_{i,j}f''$ are martingale differences as well. For $K$ big we can have
the ($L^2$) norm of $f''$ small while $f'$ is bounded. 
Without loss of generality we thus can suppose that $f$ is bounded (and from now on, we shall do so).
\medskip

For a given positive integer $v$ and positive integers $ n$,
define
$$
  F_{i,v} = \frac1{\sqrt v} \sum_{j=1}^v U_{i,j} f,\,\,\,\, i\in\Bbb Z, \,\,\,\,\text{and}\,\,\,\,
  X_{n,i} =  X_{n,i,v} = \frac1{\sqrt n} F_{i,v},\,\,\,\, i=1,\dots,n.
$$ 
Clearly, $X_{n,i}$ are martingale differences for the filtration $(\Cal F_{i, \infty})_i$. We will prove that there exists a random variable
$\eta^2$ (not necessarily at the same probability space) such that for any sequence of $v(n) \nearrow\infty$ 
$$
  E \exp \big(it  \sum_{j=1}^n X_{n,j, v(n)}\big) \to E\exp(-\frac12t^2\eta^2)
$$ 
for all $t\in \Bbb R$ as $n\to\infty$.
\medskip

The proof will use several lemmas. Their proofs will be given later.

\proclaim{Lemma 4} The sequences  $(X_{n,i,v})_i$ satisfy assumptions (i), (ii) of  Theorem 2 with $k_n=n$ uniformly for all $v\geq 1$.
\endproclaim

The next statement is adapted from \cite{V15}.

\proclaim{Lemma 5} The processes $(F_{i,v})_i$ weakly converge to a strictly stationary process $(V_i)_i$ of martingale differences defined on a 
probability space with measure $\nu_1$. The distribution of $V_1$ is a mixture of normal laws with zero means and uniformly bounded variances
and there exists a random variable $\eta_1$ such that
$$
  E V_1^2 = Ef^2, \quad E |V_1|^p <\infty,\,\,\,\,
  \frac1m \sum_{u=1}^m V_u^2 \to \eta_1^2\,\,a.s.\,(\nu_1)\,\,\,\, 
  \text{and in}\,\,\,\,L^p(\nu_1) \tag{2}
$$
for all $1\leq p<\infty$. 
\endproclaim

\proclaim{Lemma 6} There exist integers $v(m)$, $m\geq 1$, such that for any $t\in \Bbb R$, $m\to\infty$, and uniformly for all $v\geq v(m)$ 
$$
  E \exp\big( it \frac1{\sqrt m} \sum_{j=1}^m F_{j,v}\big) \to E \exp ( -\frac12 t^2\eta_1^2). \tag3
$$
\endproclaim

\medskip

Similarly as we defined the functions $F_{i,v}$ we define
$$
  G_{u,j} =  \frac1{\sqrt u} \sum_{i=1}^u U_{i,j} f
$$ 
and by the same proof as in Lemma 5 we get that there exists a probability measure $\nu_2$ on $\Bbb R^\Bbb Z$ such that for coordinate
projections $W_j$ from $\Bbb R^\Bbb Z$ to $\Bbb R$,
the processes $(G_{u,j})_j$ converge in distribution to a process $(W_j)_j$ and
$$
  \frac1m \sum_{v=1}^m W_v^2 \to \eta_2^2\,\,\,\text{a.s.}\,(\nu_2)\,\,\,\,\text{and in}\,\,\,\,L^p(\nu_2),\,\,p\geq 1
$$
where $\eta_2^2$ has all moments finite. 

\medskip

By $\Cal I_1$ let us denote the $\sigma$-algebra of sets $A\in \Cal A$ for which $T_{1,0}^{-1}A =A$. 
Suppose that $v$ is fixed. 
By Birkhoff's ergodic theorem there exists an integrable function
$$
  \eta_{v,1}^2 = E(F_{1,v}^2 \,|\, \Cal I_1) = \lim_{n\to\infty} \frac1n \sum_{i=1}^n F_{i,v}^2\,\,\,\text{a.s.}\,\,(\mu).
$$

\proclaim{Lemma 7} For $v\to\infty$, 
$$
  \eta_{v,1}^2 \overset \Cal D \to{\longrightarrow} \eta_{2}^2 
$$
where $\overset \Cal D \to{\longrightarrow}$ denotes convergence in distribution.
\endproclaim

The next statement follows from \cite{B68, Theorem 5.3}.

Let $X, X_n$, $n\geq 1$, be random variables, $p\geq 1$.
$$
  \text{If}\,\,\,\,  X_n \overset \Cal D \to{\longrightarrow} X \,\,\,\, \text{then}\,\,\,\,  E |X|^p \leq  \liminf_{n\to\infty} E |X_n|^p. \tag4
$$

Define 
$$
  Y_{n,v}^2 = \frac1n \sum_{i=1}^n F_{i,v}^2 = \sum_{i=1}^n X_{n,i,v}^2;
$$
in the proof of Lemma 6 it is shown that there exist $v(n)$, $n=1, 2, \dots$, such that for any sequence of $v_n\geq v(n)$,
$$
   Y_{n,v_n}^2 = \frac1n \sum_{i=1}^n F_{i,v_n}^2 \overset \Cal D \to{\longrightarrow} \eta_1^2 \tag{5}
$$
for $n\to\infty$. By Lemma 7, for every $n$,
$$
  E( Y_{n,v}^2\, |\, \Cal I_1) = E(F_{1,v}^2\, |\, \Cal I_1) = \eta_{v,1}^2  \overset \Cal D \to{\longrightarrow} \eta_2^2\,\,\,\text{for}\,\,\,v\to\infty.  
$$

From now on we denote $Y_n^2 =  Y_{n,v_n}^2$. Without loss of generality we can suppose $v(n) \to\infty$. 

\proclaim{Lemma 8} $Y_n^2 -  E(Y_n^2\, |\, \Cal I_1) \to 0$ in $L^1$.
\endproclaim

 From Lemma 8, Lemma 7, and (5) it follows:

\proclaim{Lemma 9} $\eta_1^2$ and $\eta_2^2$ are equally distributed.
\endproclaim

Now, we finish the proof of Theorem 1 for $d=2$. Recall that 
$$
  Y_{n,v}^2 = \frac1n \sum_{i=1}^n F_{i,v}^2 = \sum_{i=1}^n X_{n,i,v}^2,
$$
and there exist $v(n)\to \infty$ such that for any sequence of $v_n\geq v(n)$,
$$
   Y_n^2 = Y_{n,v_n}^2 = \frac1n \sum_{i=1}^n F_{i,v_n}^2 \overset \Cal D \to{\longrightarrow} \eta_1^2.
$$
Denote $\bar Y_n^2 =  E\big( Y_n^2 \,|\, \Cal I_1 \big) $. We have, for $N>n$ and $m=[N/n]$ (the integer part of $N/n$)
$$\multline
  \frac1N \sum_{i=1}^N F_{i,v}^2 = \frac1N \sum_{j=0}^{m-1} \sum_{i=1}^n F_{jn+i,v}^2 + \frac1N \sum_{i=mn+1}^N F_{i,v}^2 = \\
    \frac{mn}{N}  \bar Y_n^2 + \frac{mn}{N} \frac1m \sum_{j=0}^{m-1} \Big(\frac1n \sum_{i=1}^n F_{jn+i,v}^2  -  \bar Y_n^2\Big)  +
   \frac1N \sum_{i=mn+1}^N F_{i,v}^2.
  \endmultline
$$ 
By Lemma 8 and stationarity, $ \|(1/n)  \sum_{i=1}^n F_{jn+i,v}^2  -  \bar Y_n^2 \|_1 \to 0$ uniformly in $j$, and for $N/n \to\infty$ the last term goes 
to zero in $L^1$
as well. This proves that for any sequence of $v=v_n\geq v(n)$ and any sequence of $N_n$ with $N_n/n\to \infty$ there exist $\Cal I_1$-measurable
random variables $\eta(n)^2= \bar Y_n^2$ such that 
$$
  \big\| \eta(n)^2 - \frac1N \sum_{i=1}^N F_{i,v}^2\big\|_1 \to 0\,\,\,\,\text{and}\,\,\,\,
  \eta(n)^2 \overset \Cal D \to{\longrightarrow} \eta_1^2
$$ 
(the second statement follows from Lemma 7 and Lemma 9). \newline
The conditions (i)-(iv) and (vii) of  Proposition 3 are thus satisfied with $\eta^2 = \eta_1^2$. The random variables 
$ (1/\sqrt{n_1 n_2}) \sum_{i_1=1}^{n_1} \sum_{i_2=1}^{n_2} f\circ T_{(i_1, i_2)}$, $n_1,n_2 \to\infty$,
thus weakly converge to a law with characteristic function $\varphi(t) = \int \exp(-\eta_1^2 t^2/2)\,d\nu_1$ where the measure $\nu_1$ was 
defined in Lemma 5.

The second statement of Theorem 1 follows from (5). \newline
This finishes the proof for $d=2$.
\medskip

\centerline{\it Proofs of Lemmas 4 -  8 ($d=2$)}
\medskip

\demo{Proof of Lemma 4} The proof is well known to follow from stationarity. For reader's convenience we recall it here.

(i) For any $\epsilon >0$, uniformly for all  integers $v\geq 1$ we have
$$\multline
  \mu(\max_{1\leq i\leq n} |X_{n,i}| > \epsilon) \leq \sum_{i=1}^n \mu(|X_{n,i}| > \epsilon) =
  n \mu( | F_{0,v}| > \epsilon \sqrt n ) \leq \\
  \leq \frac1{\epsilon^2} E\Big( \Big(\frac1{\sqrt{v}} \sum_{j=1}^v U_{0,j} f\Big)^2 1_{|\sum_{j=1}^v U_{0,j} f|\geq 
  \epsilon \sqrt{nv}}\Big) \to 0
  \endmultline
$$
as $n\to \infty$. To see that the convergence is uniform for all $v$, notice that $U_{0,j}f$, $j\in\Bbb Z$, are martingale differences,
and hence by McLeish's CLT, $(1/\sqrt v) \sum_{j=1}^v U_{0,j} f$ are uniformly stochastically bounded. This proves (i).
\medskip

To see (ii) we note
$$
  \Big(\max_{1\leq i\leq n} |X_{n,i}|\Big)^2 \leq \sum_{i=1}^n X_{n,i}^2 = 
  \frac1n \sum_{i=1}^n \Big(\frac1{\sqrt v} \sum_{j=1}^v U_{i,j} f\Big)^2
$$
which implies $E\Big(\max_{1\leq i\leq n} |X_{n,i}|\Big)^2 \leq  E f^2 \leq 1$.
\enddemo \qed

\demo{Proof of Lemma 5}
For a finite set $J\subset \Bbb Z$ and for $a = (a_i;  i\in J)\in \Bbb R^J$, consider the sums
$$
  \sum_{i\in J} a_i \sum_{j=1}^v U_{i,j} f,\,\,\,\, v\to\infty.
$$
Without loss of generality (cf\. \cite{V89}) we can suppose that for the transformation $T_{0, 1}$ there exists an ergodic 
decomposition of $\mu$ into ergodic components $m_\omega$; each $m_\omega$ is a probability measure invariant and ergodic for $T_{0, 1}$.
There exists a measure $\tau$ on $(\Omega, \Cal A)$ such that for $A\in\Cal A$, $\mu(A) = \int  m_\omega(A)\, \tau(d\omega)$.
The random variables $\sum_{i\in J} a_i U_{i,j} f$, $j=1,2, \dots$, are strictly stationary martingale differences and
by  Birkhoff's ergodic theorem, 
$$
  \frac1v \sum_{j=1}^v \Big(\sum_{i\in J} a_i U_{i,j} f\Big)^2 \to \eta(a)^2\,\,\,\,\text{a.s.}\,\,(\mu)
$$
for some integrable $T_{0,1}$-invariant function $\eta(a)^2$. For almost every ergodic component $m_\omega$, $\eta(a)^2$ is a.s\. a constant
equal to $\int  \big(\sum_{i\in J} a_i U_{i,0} f\big)^2 \,dm_\omega$. By McLeish's CLT (cf\. Proposition 3) the random variables 
$ (1/\sqrt v) \sum_{j=1}^v \Big(\sum_{i\in J} a_i U_{i,j} f\Big)$ converge in distribution to $\Cal N(0, \eta(a)^2)$ (in $(\Omega, \Cal A, m_\omega)$).
By the Cramer-Wold device, for an ergodic component $m_\omega$ given and $v\to\infty$, the distributions of the 
random vectors $(F_{i,v} ; i\in J)$ thus  weakly converge to a  multidimensional normal law $\nu_\omega = \nu_{J, \omega}$.

For the measure $\mu$ we deduce that  $(F_{i,v} ; i\in J)$ weakly  converge 
to a mixture of multidimensional normal laws. 

To show this, denote $F_{J,v} = (F_{i,v} ; i\in J)$. We have shown that there are
measures $\nu_\omega$ on $\Bbb R^J$ such that for a bounded and continuous function $g$ on $\Bbb R^J$, $\int_\Omega g\circ F_{J,v}\, d\,m_\omega
\to \int_{\Bbb R^J} g\, d\nu_\omega$. We thus have  $\int_\Omega g\circ F_{J,v}\, d\mu = \int_\Omega \int_\Omega g\circ F_{J,v}\, d\,m_\omega \, 
\tau(d\omega) \to \int_\Omega  \int_{\Bbb R^J} g\, d \nu_\omega \tau(d\omega)$, and hence the measure $\nu$ on $\Bbb R^J$ defined by 
$\nu(A) =  \int_\Omega  \nu_\omega(A) \, \tau(d\omega)$ is the weak limit of $\mu\circ F_{J,v}^{-1}$ for $v\to\infty$.

This way we get a projective system of probability measures on $\Bbb R^\Bbb Z$ and following 
Kolmogorov's theorem there exists a measure $\nu_1$ on $\Bbb R^\Bbb Z$ such that for coordinate projections  $V_u$ from
 $\Bbb R^\Bbb Z$ to $\Bbb R$ and for any finite $J\subset \Bbb Z$, the vectors 
$(F_{i,v})_{i\in J}$ converge in distribution to $(V_u)_{u\in J}$ as $v\to\infty$. By strict stationarity of the sequences $(F_{u,v})_u$, $v\geq 1$,  
the measure $\nu_1$ is shift-invariant hence the process $(V_u)_u$ is strictly stationary. 

In the construction above, for $J=\{1\}$, $\nu_\omega = \Cal N(0, \int f^2 \, d\,m_\omega)$ with $ \int f^2 \, d\,m_\omega \leq \|f\|_\infty^2 <\infty$. 
Therefore  $E V_1=0$, $ E |V_1|^2 = \int \int |f|^2 \, d\,m_\omega \tau(d\omega) = \|f\|_2^2$, and for
 $1\leq p<\infty$,  $ E |V_1|^p = \int \int |x|^p \, \nu_\omega(dx) \tau(d\omega) <\infty$.

By Birkhoff's ergodic theorem there exists a random variable $\eta_1^2$ such that  for all $p\geq 1$, $\eta_1^2\in L^p(\nu_1)$ and
$(1/m) \sum_{u=1}^m V_u^2 \to \eta_1^2\,\,$a.s.$\,(\nu_1)$ and in $ L^p(\nu_1)$.
\medskip

Finally, we show that $V_i$ are martingale differences. 
It is sufficient to prove that for any bounded measurable function $g$ on $\Bbb R^n$ we have $\int V_0 g(V_{-n},\dots,V_{-1}) \,d\nu_1 = 0$.
Let us denote by $(\Cal F_k)_k$ the natural filtration of the process $(V_j)_j$.
Because $E F_{0,v}^2 = \int V_0^2\, d\nu_1$, $v\geq 1$, and the vectors $(F_{-n,v},\dots, F_{0,v})$ converge in distribution to 
$(V_{-n},\dots, V_0)$ for $v\to\infty$, $ F_{0,v}g(F_{-n,v},\dots, F_{-1,v})$ are uniformly integrable and converge in law to 
$ V_0 g(V_{-n},\dots,V_{-1})$. Therefore, $0 = E \big[ F_{0,v}g(F_{-n,v},\dots, F_{-1,v})\big] \to \int V_0 g(V_{-n},\dots,V_{-1}) \,d\nu_1$ for
 $v\to \infty$.

\enddemo
\qed

\demo{Proof of Lemma 6}

Recall that $V_j$ are stationary martingale differences in $L^2$.
 By McLeish's CLT (cf\. Proposition 3) the partial sums $(1/\sqrt m) \sum_{j=1}^m V_j$ thus converge in distribution to a law $\Cal L$ with characteristic
function $ E \exp ( -\frac12 t^2\eta_1^2)$. We can choose $m$ large enough so that the distribution of $(1/\sqrt m) \sum_{j=1}^m V_j$ is
sufficiently close to $\Cal L$ and then $v(m)$ so big that for $v\geq v(m)$ the law of $(F_{1,v},\dots, F_{m,v})$ is sufficiently close to the law
of $(V_1,\dots, V_m)$.

\enddemo
\qed

Remark that the lemma can also be proved using Theorem 2: \newline
 By Lemma 4, for $X_{m,i,v} = F_{i,v}/\sqrt m$ the assumptions (i), (ii) are satisfied (we put $k_m=m$), uniformly for all $v$.
By Lemma 5 the random vectors $(F_{1,v},\dots,F_{m,v})$ converge in law to $(V_1, \dots, V_m)$, as $v\to\infty$.
For any sequence $\ell(m)\to \infty$ with $1\leq \ell(m)\leq m$ there exist integers $v(m)$ such that
$$
  \frac1{\ell(m)} \sum_{i=1}^{\ell(m)}  F_{i,v}^2 = \eta_{\ell(m),1}^2  \to  \eta_1^2\,\,\,\,\text{in distribution, as}\,\,\,\,m\to\infty,
$$
uniformly for all $v\geq v(m)$.
 We have $E(T_m(t) \,|\, \Cal F_{\ell(m),\infty}) =  \prod_{j=1}^{\ell(m)} (1 + itX_{m,j})$. By (i), (ii) we deduce that  we can choose
$\ell(m)$ growing slowly enough so that for every $t\in \Bbb R$, $E(T_m(t) \,|\, \Cal F_{\ell(m),\infty})  \to 1$ in $L^1$ as $m\to\infty$. For 
$\eta_m^2 =  \eta_{\ell(m),1}^2$ the conditions (iii), (iv), (v) of  Theorem 2 are thus satisfied.
\qed

\demo{Proof of Lemma 7} Let $t\in \Bbb R$.
By McLeish's theorem (cf\. \cite{HaHe, Theorem 3.2} or Theorem 2), for $v$ fixed and $n\to\infty$ 
$$
  E \exp\Big(it\frac1{\sqrt n} \sum_{l=1}^n F_{l,v}\Big) \to E \exp( -\frac12 t^2\eta_{v,1}^2). \tag{3a}
$$
By definition we have
$$
 \frac1{\sqrt n}  \sum_{i=1}^n  F_{i,v} =  \frac1{\sqrt v}  \sum_{j=1}^v  G_{n,j}
$$
and by Lemma 6 applied to the functions $G_{n,j}$, for every $v$ there exists an $n(v)$ such that for $v\to \infty$
$$
  E \exp \Big(it\frac1{\sqrt v} \sum_{j=1}^v G_{n,j}\Big) \to E \exp( -\frac12 t^2\eta_{2}^2)
$$
uniformly  for $n\geq n(v)$.
For a given $\epsilon>0$ we thus can find $v_0$ big enough so that for all $v\geq v_0$ and all $n\geq n(v)$, 
$$
  |  E \exp \Big(it\frac1{\sqrt v} \sum_{j=1}^v G_{n,j}\Big) - E \exp( -\frac12 t^2\eta_{2}^2)| <\epsilon ;
$$ 
keeping $v$ fixed, by (3a) there is $n\geq n(v)$ large enough so that
$$
  |  E \exp\Big(it\frac1{\sqrt n} \sum_{j=1}^n F_{j,v}\Big) - E \exp( -\frac12 t^2\eta_{v,1}^2)| < \epsilon.
$$ 
Therefore, $  |E \exp( -\frac12 t^2\eta_{2}^2) -  E \exp( -\frac12 t^2\eta_{v,1}^2)| < 2\epsilon$.
By properties of the  Laplace transformation (cf\. \cite{F, Chapter XIII}), the convergence 
$  E \exp( -\frac12 t^2\eta_{v,1}^2) \to E \exp( -\frac12 t^2\eta_{2}^2) $,  $v\to\infty$, $t\in \Bbb R$, implies
$\eta_{v,1}^2 \overset \Cal D \to{\longrightarrow} \eta_{2}^2 $.
\enddemo
\qed

\demo{Proof of Lemma 8}
For $k\geq 0$ define
$$
  \varphi_k(x) = \cases  x\,\,\,&\text{if}\,\,\, |x| \leq k, \\
                                      k\cdot sign(x) &\text{if}\,\,\, |x|>k.
                          \endcases 
$$
By Jensen's inequality taken conditionally we have, for $Z= \varphi_k(Y_n^2)$, $E(Z^2\,|\,\Cal I_1) \geq [E(Z\,|\,\Cal I_1)]^2$ a.s\. hence for every 
$k, n\geq 1$,
$$
  E\Big[\varphi_k(Y_n^2) \Big]^2 \geq E\Big[ E\big( \varphi_k(Y_n^2)\,|\,\Cal I_1\big) \Big]^2.  \tag{6}
$$
By (5),
$
   E\Big[\varphi_k(Y_n^2) \Big]^2 \underset n\rightarrow \infty \to{\longrightarrow} E\Big[\varphi_k(\eta_1^2) \Big]^2
$
for every $k$. We thus have
$$
  \lim_{k\to\infty} \lim_{n\to \infty} E\Big[\varphi_k(Y_n^2) \Big]^2  =  E \eta_1^4.   \tag{7}
$$
By concavity of $\varphi_k$ (on $[0, \infty)$), $(1/n) \sum_{i=1}^n \varphi_k(F_{i,v_n}^2) \leq \varphi_k\big( (1/n)  \sum_{i=1}^n F_{i,v_n}^2 \big) = \varphi_k(Y_n^2)$; we deduce
that for every $k$,
$$
   E\big(\varphi_k(F_{1,v_n}^2)  \,|\, \Cal I_1\big) = E\Big(\frac1n \sum_{i=1}^n \varphi_k(F_{i,v_n}^2)  \,|\, \Cal I_1\Big) \leq 
   E\big(\varphi_k(Y_n^2) \,|\, \Cal I_1\big).  \tag{8}
$$

Recall that by Lemma 5 and Lemma 7, for $n\to\infty$,
$$
 F_{1,v_n} \overset \Cal D \to{\longrightarrow} V_1,   \,\,\,\, E  F_{1,v_n}^2 =E f^2 = E V_1^2,\,\,\,\,\, 
  E(F_{1,v_n}^2\, |\, \Cal I_1)  \overset \Cal D \to{\longrightarrow} \eta_2^2.
$$
By uniform integrability of $  F_{1,v_n}^2$, for every $\epsilon>0$ there exist a $k\geq 1$ and $n(k)$ such that for all $n\geq n(k)$,
$
  E \big| F_{1,v_n}^2 - \varphi_k(F_{1,v_n}^2)\big| < \epsilon.
$
By contractiveness of conditional expectation we get $ E \big|   E(F_{1,v_n}^2\, |\, \Cal I_1)  - E( \varphi_k(F_{1,v_n}^2)  \,|\, \Cal I_1) \big| <\epsilon$.
For a given $\delta>0$ we therefore can choose $\epsilon>0$ small enough so that
$$
  \mu\big\{ |  E(F_{1,v_n}^2\, |\, \Cal I_1)  - E( \varphi_k(F_{1,v_n}^2)  \,|\, \Cal I_1)| >\delta \big\} <\delta. \tag{9}
$$
From $\big[ E(F_{1,v_n}^2  \,|\, \Cal I_1)\big]^2  \overset \Cal D \to{\longrightarrow} \eta_2^4$ and (4) it follows
$$
   \liminf_{n\to\infty} E[E(F_{1,v_n}^2\, |\, \Cal I_1)]^2 \geq E \eta_2^4.
$$
By (2), $ E \eta_2^4 <\infty$ (we can prove (2) for $\eta_2$ in the same way as for $\eta_1$).
By (9) we have $\lim_{k\to\infty,\, n\geq n(k)} E( \varphi_k(F_{1,v_n}^2)  \,|\, \Cal I_1)^2  \overset \Cal D \to{\longrightarrow}  \eta_2^4$
hence using (4) again we deduce that
for every $\epsilon>0$ there are $k(\epsilon)$ and $n(\epsilon, k)$ such that for all 
$k\geq k(\epsilon), n\geq n(\epsilon, k)$, 
$$
    E\Big[E\big(\varphi_k(F_{1,v_n}^2)  \,|\, \Cal I_1\big) \Big]^2 \geq  E \eta_2^4 -\epsilon.
$$
Therefore,
$$
  \lim_{k\to\infty} \liminf_{n\to\infty}  E\Big[E\big(\varphi_k(F_{1,v_n}^2)  \,|\, \Cal I_1\big) \Big]^2 \geq  E \eta_2^4. \tag{10}
$$
From (7), (6), (8), (10)  we deduce that for every $\epsilon>0$ there are $k(\epsilon)$ and $n(k,\epsilon)$ such that for $k\geq k(\epsilon)$  and 
$n\geq n(k,\epsilon)$ 
$$\multline
  E \eta_1^4 +\epsilon \geq  E\big[\varphi_k(Y_n^2) \big]^2 \geq  E\Big[ E\big( \varphi_k(Y_n^2) \,|\, \Cal I_1 \big) \Big]^2 \geq \\
  E\Big[ E( \varphi_k(F_{1,v_n}^2) \,|\, \Cal I_1 \big) \Big]^2 \geq E \eta_2^4 -\epsilon. \endmultline \tag{11}
$$
Therefore, $E \eta_1^4  \geq  E \eta_2^4$ and by symmetry we get
$$
  E \eta_1^4  =  E \eta_2^4.
$$
Using (11) we deduce
$$  
   \lim_{k\to\infty} \limsup_{n\to \infty} \Big\{ E\big[\varphi_k(Y_n^2) \big]^2 - E\big[ E\big( \varphi_k(Y_n^2)\,|\,\Cal I_1\big) \big]^2 \Big\} = 0.
$$
This implies
$$
   \lim_{k\to\infty} \limsup_{n\to \infty} E\big[\varphi_k(Y_n^2) -  E\big( \varphi_k(Y_n^2) \,|\, \Cal I_1 \big) \big]^2 = 0,
$$
hence
$$
   \lim_{k\to\infty} \limsup_{n\to \infty} \big\|\varphi_k(Y_n^2) -  E\big( \varphi_k(Y_n^2) \,|\, \Cal I_1 \big) \big\|_1 = 0.
$$
Because $E Y_n^2 = E f^2 = E \eta_1^2$ and $Y_n^2 \overset \Cal D \to{\longrightarrow} \eta_1^2$, $Y_n^2$ are uniformly integrable, 
hence 
$$
   \lim_{k\to\infty} \limsup_{n\to \infty}  \| Y_n^2 - \varphi_k(Y_n^2) \|_1 = 0.
$$
By triangular inequality and contractiveness of conditional expectation we get
$$\gather
   \lim_{k\to\infty} \limsup_{n\to \infty}  \big\| Y_n^2 -  E\big( \varphi_k(Y_n^2) \,|\, \Cal I_1 \big) \big\|_1 = 0, \\
   \lim_{k\to\infty} \limsup_{n\to \infty}  \big\|  E\big( Y_n^2 \,|\, \Cal I_1 \big) -  E\big( \varphi_k(Y_n^2) \,|\, \Cal I_1 \big)\big\|_1 = 0,
  \endgather
$$ 
therefore
$$
   \lim_{n\to \infty} \big\|  Y_n^2 -  E\big( Y_n^2 \,|\, \Cal I_1 \big) )\big\|_1 = 0;
$$
this finishes the proof of Lemma 8.

\enddemo
\qed

\bigskip

\centerline{\it The case of $d>2$}
\bigskip

Now we give an idea of a  proof of the theorem for $d>2$. Suppose that $d> 2$ and for dimension $d-1$, the theorem has been proved.
 Using the same reasoning  as before we can restrict ourselves to the case of $f$ bounded. 

For each $k\in \{1,2,\dots,d\}$, let us denote by $\underline{i^k}$ the vectors  $\underline{i} =(i_1,i_2,\dots,i_d) \in \Bbb Z^d$ such that $i_k=0$. 
To $\underline{i}\in \Bbb Z^{d-1}$ we thus can associate a $\underline{i^k}\in \Bbb Z^d$ and and having a 
$\Bbb Z^d$ action we thus get a $ \Bbb Z^{d-1}$ action  by associating $\underline{i} \mapsto T_{\underline{i^k}}$. By the induction hypothesis 
(cf\. the proof for $d-1$), by boundedness of $f$,
the sums $(1/\sqrt{|\underline{n^k}}|) \sum_{\underline{1^k}\leq \underline{j^k}\leq \underline{n^k}} f\circ T_{\underline{j^k}}$ converge 
in distribution to a law with characteristic function $ E\big(-\frac12 t^2 \eta^2\big)$ with $\eta \in L^p$ for all $1\leq p<\infty$.

We will prove a convergence of $(1/\sqrt{\underline{n}}) \sum_{\underline{1}\leq \underline{j}\leq \underline{n}} f\circ T_{\underline{j}}$ 
($\underline{n} \in \Bbb Z^d$, $\underline{n}\to \underline{\infty}$)
to a  limit law with characteristic function $\varphi(t) = E \exp(-\bar \eta^2 t^2/2)$ where $\bar \eta^2 \in L^p$ for all $1\leq p<\infty$.
By $(i, \underline{v})$ we denote an element of $\Bbb Z^{d}$ where $i\in \Bbb Z$ and $\underline{v} \in \Bbb Z^{d-1}$. We define 
$$
  F_{i, \underline{v}} = \frac1{\sqrt{| \underline{v}|}} \sum_{ \underline{j}=  \underline{1}}^{ \underline{v}} U_{j, \underline{v}} f,
  \,\,\,\,  X_{n,i} = X_{n,i,\underline{v}} = \frac1{\sqrt{| \underline{n}|}}  F_{i, \underline{v}}, \,\,\,\,i=1,\dots,n.
$$
For $\underline{v}$ fixed, $ F_{i, \underline{v}}$ are martingale differences for the filtration $(\Cal F_i^{(1)})_i$. Similarly as in the preceding case
we will prove that there exists a random variable $\eta^2$ such that for every $t\in \Bbb R$ and any sequence $\underline{v}(n) \to \underline{\infty}$,
$$
  E \exp \big(it \sum_{j=1}^n X_{n,j,\underline{v}(n)}\big) \to E\big(-\frac12 t^2 \eta^2\big)
$$
as $n\to\infty$.

Lemma 4 remains valid for $ X_{n,i,\underline{v}} $ with the same proof:

\proclaim{Lemma 4a} The sequences  $(X_{n,i,\underline{v}})_i$ satisfy assumptions (i), (ii) of  Theorem 2 with $k_n=n$ uniformly for all 
$\underline{v}\geq \underline{1}$.
\endproclaim

For $u\in \Bbb Z$ and $\underline{j} \in \Bbb Z^{d-1}$ we define
$$
  G_{u,\underline{j}} = \frac1{\sqrt u} \sum_{i=1}^u U_{i, \underline{j}} f.
$$

\proclaim{Lemma 5a} For $\underline{v}\to \underline{\infty}$ ($\underline{v}\in \Bbb Z^{d-1}$) the sequences $( F_{u, \underline{v}})_u$ converge 
in distribution to a sequence of martingale differences $(V_u)_{u\in \Bbb Z}$ where $EV_1^2 = Ef^2$ and $V_1\in L^p$ 
for all $1\leq p<\infty$. We have $(1/m) \sum_{u=1}^m V_u^2 \to \eta_1^2$ a.s\. and in $L^p$, $1\leq p<\infty$. 

The random fields $(  G_{u,\underline{j}})_{\underline{j}}$ converge in distribution to a stationary random field $(W_{\underline{j}})_{\underline{j}}$
defined on a space with a probability measure $\nu_2$. For every $1\leq p<\infty$, $W_{\underline{j}} \in L^p(\nu_2)$,  $E W_1^2 = Ef^2 = 
E G_{u,\underline{j}}^2$, and $(W_{\underline{j}})_{\underline{j}}$ is a multiple field of martingale differences; 
$W_{\underline{j}}\in L^p$ for all  $1\leq p<\infty$.
\endproclaim

\demo{Proof} The first part of the Lemma (existence of the sequence $(V_u)_u$) can be proved in the same way as Lemma 5, only instead of
 McLeish's theorem we use the CLT for $(d-1)$-dimensional random fields from our assumptions. By assumptions, as mentioned above, the law of $V_u$
has characteristic function $ E\big(-\frac12 t^2 \eta^2\big)$ with $\eta \in L^p$ for all $1\leq p<\infty$ hence $V_u\in L^p$.

To prove the second part of the Lemma we can in the same way as before, using  McLeish's theorem, show that the random fields 
$( G_{u,\underline{j}})_{\underline{j}}$ converge in distribution to a stationary random field $(W_{\underline{j}})_{\underline{j}}$ of coordinate 
projections of $\Bbb R^{\Bbb Z^{d-1}}$ (onto $\Bbb R$) with the product 
$\sigma$-algebra of Borel sets and a probability measure $\nu_2$ invariant with respect to the shifts;  $W_{\underline{j}} \in L^p(\nu_2)$ for every 
$1\leq p<\infty$, $E W_1^2 = Ef^2 = E G_{u,\underline{j}}^2$ (for every $u, \underline{j}$).

Let us show that $(W_{\underline{j}})_{\underline{j}}$ is a  multiple field of martingale differences. The $\sigma$-algebras $ \Cal F_i^{(q)}$ 
generated by the random variables $W_{\underline{j}}$ with $j_q\leq i$ are a multiple filtration. By definition, $W_{\underline{j}}$ is 
$ \Cal F_{j_q}^{(q)}$-measurable for all $1\leq q\leq d$ and using the same argument as in the proof of Lemma 5 we can show that 
$E(W_{\underline{j}}\,|\,  \Cal F_{j_q-1}^{(q)}) =0$.

In the same way as in the proof of Lemma 5 we can show that  $EW_{\underline{j}}^2 =Ef^2$ and $W_{\underline{j}}\in L^p$ for all  $1\leq p<\infty$. 

\enddemo
\qed

By Theorem 1 for dimension $d-1$, for the random field $( W_{\underline{j}})_{\underline{j}}$ there exists a measurable function $\eta_2^2$ such that 
for $\underline{n} \to \underline{\infty}$, $(1/\sqrt{\underline{n}}) \sum_{\underline{1}\leq \underline{j}\leq \underline{n}} f\circ T_{\underline{j}}$
converge in distribution to a limit law with characteristic function $\varphi(t) = \int \exp(-\eta_2^2 t^2/2)\,d\nu_2$.

In the same way as for $d=2$ we can prove the following lemma.

\proclaim{Lemma 6a}
There exist $\underline{v}(m)\in \Bbb Z ^{d-1}$, $m\geq 1$, such that for $m\to\infty$ and uniformly for all $\underline{v}\geq \underline{v}(m)$ 
$$
  E \exp\big( it \frac1{\sqrt m} \sum_{l=1}^m F_{l,\underline{v}}\big) \to E \exp ( -\frac12 t^2\eta_1^2), \,\,\,\,t\in \Bbb R.
$$
There exist $u(\underline{v})\geq 1$ and $\underline{v}\in \Bbb Z ^{d-1}$, $\underline{v}\geq \underline{1},$ such that for 
$\underline{v}\to \underline{\infty}$ and uniformly for all $u\geq u(\underline{v})$ 
$$
  E \exp \Big(it \frac1{\sqrt{\underline{v}}} \sum_{\underline{1} \leq \underline{j}\leq \underline{v}} U_{\underline{j}} G_{u,\underline{j}}\Big)
  \to  E \exp ( -\frac12 t^2\eta_2^2), \,\,\,\,t\in \Bbb R.
$$
\endproclaim

By $\Cal I_1$ let us denote the $\sigma$-algebra of sets $A\in \Cal A$ for which $T_{\underline{e_1}}^{-1}A =A$. 
Suppose that $\underline{v}\in \Bbb Z^{d-1}$ is fixed. 
By Birkhoff's ergodic theorem there exists an integrable function
$$
  \eta_{\underline{v},1}^2 = E(F_{1,\underline{v}}^2 \,|\, \Cal I_1) = \lim_{n\to\infty} \frac1n \sum_{i=1}^n F_{i,\underline{v}}^2.
$$

\proclaim{Lemma 7a} For $\underline{v}\to\underline{\infty}$, 
$$
  \eta_{\underline{v},1}^2 \overset \Cal D \to{\longrightarrow} \eta_{2}^2.
$$
\endproclaim

Lemma 7a can be proved in the same way as Lemma 7.
\medskip

As in the case of $d=2$ we define 
$$
  Y_{n,\underline{v}}^2 = \frac1n \sum_{i=1}^n F_{i,\underline{v}}^2 = \sum_{i=1}^n X_{n,i,\underline{v}}^2;
$$
there exist $\underline{v}(n)$, $n=1, 2, \dots$, such that for any sequence of $\underline{v}_n\geq \underline{v}(n)$,
$$
   Y_{n,\underline{v}_n}^2 = \frac1n \sum_{i=1}^n F_{i,\underline{v}_n}^2 \overset \Cal D \to{\longrightarrow} \eta_1^2 \tag{5a}
$$
for $n\to\infty$. By Lemma 7, for every $n$,
$$
  E( Y_{n,\underline{v}}^2\, |\, \Cal I_1) = E(F_{1,\underline{v}}^2\, |\, \Cal I_1) = \eta_{\underline{v},1}^2  \overset \Cal D \to{\longrightarrow} 
  \eta_2^2\,\,\,\text{for}\,\,\,\underline{v}\to\underline{\infty}.  
$$

From now on we denote $Y_n^2 =  Y_{n,\underline{v}_n}^2$. Without loss of generality we can suppose $\underline{v}(n) \to\underline{\infty}$. 
Now, Lemma 8 can be stated and proved as for $d=2$ and then the rest of the proof of Theorem 1 can be done as before.
Notice that the assumption of $f$ bounded leads in the CLT to a limit distribution with characteristic function $ E \exp ( -\frac12 t^2\eta_1^2)$
where $\eta_1\in L^p$ for all $1\leq p<\infty$.

\enddemo
\qed

\subheading{3. A weak invariance principle}

By a weak invariance principle (WIP) we shall mean the WIP as defined in \cite{VWa}. For $\underline{n}= (n_1,\dots,n_d)\in \Bbb N^d$ and 
$\underline{t}= (t_1,\dots, t_d)\in [0, 1]^d$ we define $\underline{n} \cdot \underline{t} = (n_1t_1,\dots n_dt_d)$,
 $[\underline{n} \cdot \underline{t}] = ([n_1t_1],\dots, [n_dt_d])$; we denote
$$
  S_{\underline{n},\underline{t}}(f) = \sum_{\underline{0}\leq \underline{k}\leq [\underline{n} \cdot \underline{t}]} \prod_{i=1}^d
  (k_i  \wedge (n_it_i - 1) -k_i+1) U_{ \underline{k}}f, \quad s_{\underline{n},\underline{t}}(f) = S_{\underline{n},\underline{t}}(f) / 
  \sqrt{|\underline{n}|}.
$$
Notice that for $d=1$, $ S_{\underline{n},\underline{t}}(f) = \sum_{k=0}^{[nt]} U_kf + (nt-[nt]) U_{[nt]+1}f$ (here we denote $U_k = U^k$).

\proclaim{Theorem 10}  Let $f\in L^2$, be such that $(f\circ T_{\underline{i}})_{\underline{i}}$ ($\underline{i}\in \Bbb Z^d$)  is a field of martingale 
differences for a completely commuting filtration $(\Cal F_{\underline{i}})_{\underline{i}}$ ($\underline{i}\in \Bbb Z^d$). If $n_j \to\infty$, $j=1,\dots, d$, 
then there exists a random 
variable $\eta\geq 0$ and a Brownian sheet $(\Bbb B_{\underline{t}})_{\underline{t}\in [0, 1]^d}$ independent of $\eta$ such that for 
$n_1,\dots, n_d \to \infty$ $ s_{\underline{n},\underline{t}}(f)$  converge in distribution to $(\eta \Bbb B_{\underline{t}})_{\underline{t}}$.

\endproclaim

\demo{Proof} The proof can be done in the same way as in \cite{CDV}. We have to prove convergence of finite dimensional distributions
and tightness. For tightness, the proof from \cite{VWa} remains valid (note that in \cite{CDV} the same proof, adapted to reverse martingales, is used). 

The  convergence of finite dimensional distributions will be proved for $d=2$; for $d>2$ the result can be extended by induction.
As in \cite{CDV}  it is necessary and sufficient to show that for partitions $0=t_0<t_1<\dots<t_K=1$ and $0=s_0<s_1<\dots<s_K=1$ of the interval 
$[0, 1]$  and real constants $a_{k,l}$, $0\leq k, l\leq K-1$, and $n_1,n_2 \to\infty$, the sums 
$$
  \frac1{\sqrt{n_1n_2}} \sum_{k=0}^{K-1}  \sum_{l=0}^{K-1} a_{k,l} \sum_{i=[n_1t_k]}^{[n_1t_{k+1}]-1} \sum_{j=[n_2s_l]}^{[n_2s_{l+1}]-1}
  U_{i,j}f
$$
weakly converge to 
$$
  \eta \sum_{k=0}^{K-1}  \sum_{l=0}^{K-1} a_{k,l} (\Bbb B_{t_{k+1}, s_{l+1}} - \Bbb B_{t_{k+1}, s_{l}} - \Bbb B_{t_{k}, s_{l+1}}
  + \Bbb B_{t_{k}, s_{l}} )
$$
where $\Bbb B_{t,s}$ is a  Brownian sheet independent of $\eta$; $\eta\geq 0$ is defined as in the proof of Theorem 1.

Let us fix a $0\leq k\leq K-1$,
let $J\subset \Bbb Z$ be a finite set and for $\bar a = (\bar a_i;  i\in J)$ ($\bar a_i\in \Bbb R$ are constants), consider, as  in the proof of Lemma 5, the sums
$$
  \sum_{i\in J} \bar a_i  \sum_{l=0}^{K-1} a_{k,l} \sum_{j=[t_lv]}^{[t_{l+1}v]-1}  U_{i,j} f,\,\,\,\, v\to\infty.
$$
As in the proof of Lemma 5, let $m_\omega$ be invariant and ergodic probability measures for the transformation $T_{0, 1}$ (ergodic decomposition
of $\mu$).
By Birkhoff's ergodic theorem
$$\gather
   \frac1v \sum_{j=1}^v \Big(\sum_{i\in J} a_i U_{i,j} f\Big)^2 \to  \eta(\bar a)^2, \\
  \frac1v \sum_{l=0}^{K-1}  a_{k,l}^2  \sum_{j=[t_lv]}^{[t_{l+1}v]-1}  \Big( \sum_{i\in J} \bar a_i  U_{i,j} f\Big)^2 \to
   \eta(\bar a)^2 \sum_{l=0}^{K-1} (t_{l+1}-t_l) a_{k,l}^2 
\endgather
$$
(a.s\. for $\mu$ hence also for almost all  $m_\omega$).
 In the same way as in the proof of Lemma 5 we conclude that
$$
  (1/\sqrt v) \sum_{i\in J} \bar a_i  \sum_{l=0}^{K-1} a_{k,l} \sum_{j=[t_lv]}^{[t_{l+1}v]-1}  U_{i,j} f \overset \Cal D \to{\longrightarrow}
  \Cal N\big(0, \eta(\bar a)^2 \sum_{l=0}^{K-1} (t_{l+1}-t_l) a_{k,l}^2\big).
$$
In the same way as in the proof of Lemma 5 we deduce, using the Cramer-Wold 
device and Kolmogorov's theorem on projective limit, that the processes 
$$
  \Big(\frac1{\sqrt v}  \sum_{l=0}^{K-1} a_{k,l} \sum_{j=[t_lv]}^{[t_{l+1}v]-1}  U_{i,j} f\Big)_i
$$ 
converge in distribution to a process $( (\sum_{l=0}^{K-1} (t_{l+1}-t_l) a_{k,l}^2)^{1/2} V_i)_i$
where $(V_i)_i$ is the process found in the proof of Lemma 5.

In the same way as in the proof of Theorem 1 we deduce that for $n_1, n_2\to\infty$ the sums
$$
  \frac1{\sqrt{n_1n_2}}  \sum_{i=1}^{n_1} \sum_{l=0}^{K-1} a_{k,l} \sum_{j=[t_ln_2]}^{[t_{l+1}n_2]-1}  U_{i,j} f
$$
converge in distribution to the law with characteristic function $ E\exp (-\eta^2c_k^2t^2/2)$ where $\eta^2$ is the same random variable as in Theorem 1
and $ c_k^2 = \sum_{l=0}^{K-1} (t_{l+1}-t_l) a_{k,l}^2 $.

Recall that $\Cal I_1$ is the $\sigma$-algebra of $T_{1,0}$-invariant sets from $\Cal A$ and denote
$$
  Y_{k,n_1,n_2}^2 = \frac1{n_1}  \sum_{i=1}^{n_1} \frac1{n_2} \Big( \sum_{l=0}^{K-1} a_{k,l} \sum_{j=[t_ln_2]}^{[t_{l+1}n_2]-1}  U_{i,j} f\Big)^2.
$$

In the same way as in Lemma 8 we can see that for any $\epsilon>0$ there are $n_1(\epsilon), n_2(\epsilon)$ such that for 
$n_1\geq n_1(\epsilon), n_2\geq  n_2(\epsilon)$ we have $\|  Y_{k,n_1,n_2}^2 - E( Y_{k,n_1,n_2}^2 \,|\,\Cal I_1)\|_1 <\epsilon$.

From this we deduce that 
$$\multline
  \frac1{n_1n_2} \sum_{k=0}^{K-1}  \sum_{i=[n_1t_k]}^{[n_1t_{k+1}]-1} \Big( \sum_{l=0}^{K-1} a_{k,l} \sum_{j=[n_2s_l]}^{[n_2s_{l+1}]-1}
  U_{i,j}f\Big)^2 \to   \\
  \eta^2 \sum_{l=0}^{K-1}  \sum_{k=0}^{K-1} (t_{l+1}-t_l)  (s_{k+1}-s_k)    a_{k,l}^2
  \endmultline
$$
in $L^1$ for $n_1, n_2 \to \infty$, in the sense that  for any $\epsilon>0$ there are $n_1(\epsilon), n_2(\epsilon)$ such that for $n_1\geq n_1(\epsilon), 
n_2\geq n_2(\epsilon)$ the norm is smaller than $\epsilon$.

Because Lemma 4 applies to the random variables 
$$
   X_{n_1,n_2,i} = \frac1{\sqrt{n_1n_2}}  \sum_{l=0}^{K-1} a_{k,l} \sum_{j=[n_2s_l]}^{[n_2s_{l+1}]-1}  U_{i,j}f, \quad [n_1t_k]\leq i < [n_1t_{k+1}],
  \,\, 0\leq k<K,
$$
by Proposition 3
$$
  \frac1{\sqrt{n_1n_2}} \sum_{k=0}^{K-1}  \sum_{l=0}^{K-1} a_{k,l} \sum_{i=[n_1t_k]}^{[n_1t_{k+1}]-1} \sum_{j=[n_2s_l]}^{[n_2s_{l+1}]-1}
  U_{i,j}f
$$
converge in distribution to the law with characteristic function  $ E\exp (-\eta^2c^2t^2/2)$ where
$$
  c^2 =  \sum_{l=0}^{K-1}  \sum_{k=0}^{K-1} (t_{l+1}-t_l)  (s_{k+1}-s_k)    a_{k,l}^2.
$$

\enddemo
\qed

\bigskip

\subheading{4. Applications}

\centerline{\it Central limit theorems}
\bigskip

Many limit theorems have been proved using approximations by martingales. Let us denote, for $\underline{n} = (n_1,\dots,n_d)$, 
$|\underline{n}| = n_1\dots n_d$ and $S_{\underline{n}}(f) = \sum_{i_1=1}^{n_1} \dots \sum_{i_d=1}^{n_d} f\circ T_{(i_1,\dots,i_d)}$.
Using Theorem 1 we can in the same way as in \cite{Go69} prove that if $m\circ T{\underline{i}}$ are martingale differences in $L^2$
($\underline{i} \in \Bbb Z^d$), and for $n_1,\dots, n_d\to\infty$ 
$$
  \frac1{\sqrt{|\underline{n}|}} \big\| S_{\underline{n}}(f-m)\big\|_2 \to 0 \tag{12}
$$
then there exists a random variable $\eta^2$ (possibly in another probability space), $E\eta^2 =\|m\|_2^2$, such that  
$(1/\sqrt{|\underline{n}|})  S_{\underline{n}}(f)$ weakly converge to the law with characteristic function $E \exp(-\frac12 t^2\eta^2)$.

The CLT is often accompanied with a weak inviariance principle (not always, cf\. e.g\. \cite{V-Samek}, \cite{DMV}, \cite{GiV}). We will present 
conditions leading to (12) and a  weak inviariance principle at the same time.
\bigskip
\centerline{\it Weak invariance principles (WIP)}
\bigskip

A ``closer" approximation than (12) guarantees the WIP. For $d=1$ a classical assumption has been the martingale-coboundary decomposition
$f = m + g-Ug$ where $m, g\in L^2$ and $(U^ig)_i$ is a sequence of martingale differences. For $d>1$ we define the martingale-coboundary 
representation
$$
  f = \sum_{S\subset \{1,\dots,d\}}\prod_{q\in S^c} (I-U_{\underline{e_q}})g_S\,\,\,\,(f\in L^p, \,\, 1\leq p<\infty) \tag{13}
$$
where for $S\subset \{1,\dots,d\}$, $g_S \in \underset q\in S \to{\bigcap}  L^p(\Cal F_{0}^{(q)}) \ominus  L^p(\Cal F_{-1}^{(q)})$;
$\prod_{q\in \emptyset} (I-U_{\underline{e_q}})$ is defined as $I$, the identity operator. For $q'\in S$, 
$U_{\underline{e_{q'}}}^i\prod_{q\in S^c} (I-U_{\underline{e_q}})g_S$, $i\in \Bbb Z$, are martingale differences while for $q'\in S^c$, 
$\prod_{q\in S^c} (I-U_{\underline{e_q}})g_S$ are coboundaries for the transformation $T_{\underline{e_{q'}}}$. Eq\. (13) with all terms in $L^2$ implies
the WIP. For the first time, this decomposition was studied by Gordin in \cite{Go09}. Sufficient and necessary and sufficient conditions for (13)
can be found in \cite{EGi}, \cite{V17}, \cite{Gi}).

If $f = \sum_{\underline{i}\in \Bbb Z^d} P_{\underline{i}}f$ ($P_{\underline{i}}$  are the orthogonal projection operators defined in Section 1 - Introduction)
then we say that $f$ (or the random field of $f\circ  T_{\underline{i}}$) is {\it regular}.

As we will state in Theorem 11, for regular random fields, (12) and a WIP are guaranteed e.g\. by the conditions of Hannan (see \cite{VWa})
$$
   \text{(Hannan)}\quad  f = \sum_{\underline{j}\in \Bbb Z^d} P_{\underline{j}}f, \quad  \sum_{\underline{j}\in \Bbb Z^d} \|P_{\underline{j}}f\|_2 <\infty 
  \tag{14}
$$
and  Maxwell-Woodrooffe (see \cite{Gi}; the CLT was proved in \cite{PZ})
$$
   \text{(Maxwell-Woodroofe)}\quad  f = \sum_{\underline{j}\in \Bbb Z^d} P_{\underline{j}}f, \quad \sum_{\underline{u}\geq \underline{1}}
  \frac{\| E(S_{\underline{u}}\,|\,\Cal F_{\underline{1}})\|_2}{|\underline{u}|^{3/2}} <\infty. \tag{15}
$$

\proclaim{Theorem 11} If $ (T_{\underline{i}})_{\underline{i}}$ is a $\Bbb Z^d$ action, the random field $(f\circ  T_{\underline{i}})_{\underline{i}}$
is regular and one of the conditions (13), (14), (15) holds ((13) with all terms in $L^2$) then the WIP holds.
\endproclaim

\demo{Proof} In  \cite{EGi}, \cite{V17}, \cite{VWa},  \cite{Gi} it was shown that each of the conditions implies 
$$
  S_{\underline{n}, \underline{t}}(f) =  S_{\underline{n}, \underline{t}}(m) + R_{\underline{n}, \underline{t}}(f), \quad  \underline{t}\in [0, 1]^d,
$$
where $m\in L^2$,  $ U_{\underline{i}}m$ are martingale differences, and $R_{\underline{n}, \underline{t}}(f)$ converge in distribution to the zero
function on $[0, 1]^d$. The WIP thus holds for  the random field of $ U_{\underline{i}}f$ if and only if it holds for the random field of $ U_{\underline{i}}m$.     For $ U_{\underline{i}}m$ the WIP is guaranteed by  Theorem 10.

\enddemo
\qed 

\subheading{Thanks} I greatly thank the two unknown (to me) referees for detailed and careful reading of the manuscript and for many helpful comments and remarks.

\end